\documentclass[11pt]{amsart}

\usepackage[english]{babel}
\usepackage{amsmath}
\usepackage{amsfonts} 
\usepackage{amsthm}
\usepackage{times}
\usepackage{hyperref}
\usepackage{graphicx}

\newcommand{\N}{\mathbb N}

\newcommand{\R}{\mathbb R}

\newcommand{\MS}{{\mathcal S}}
\newcommand{\MN}{{\mathcal N}}
\newcommand{\MF}{{\mathcal F}}

\newcommand{\MI}{{\mathcal I}}
\newcommand{\MU}{{\mathcal U}}

\DeclareMathOperator{\interior}{int}
\DeclareMathOperator{\bd}{bd}
\DeclareMathOperator{\pos}{pos}

\DeclareMathOperator{\vol}{vol}
\DeclareMathOperator{\conv}{conv}

\newtheorem{definition}{Definition}
\newtheorem{proposition}{Proposition}
\newtheorem{theorem}{Theorem}
\newtheorem{corollary}{Corollary}
\newtheorem*{problem}{Problem}

\renewcommand{\vec}[1]{\boldsymbol{#1}}

\author{Achill Sch\"urmann}
\address{Mathematics Department, University of Magdeburg, 39106 Magdeburg, Germany}
\email{achill@math.uni-magdeburg.de}
\title{On packing spheres into containers\\(about Kepler's finite sphere packing problem)}
\date{\today}

\begin{document}

\begin{abstract}
In an Euclidean $d$-space, the container problem 
asks to pack $n$ equally sized spheres into 
a minimal dilate of a fixed container. 
If the container is a smooth convex body
and $d\geq 2$ we show that solutions to the container 
problem can not have a ``simple structure'' for large $n$.
By this we in particular find that there exist arbitrary small $r>0$,
such that packings in a smooth, $3$-dimensional convex body,
with a maximum number of spheres of radius $r$,
are necessarily not hexagonal close packings.
This contradicts Kepler's famous statement
that the cubic or hexagonal close packing 
``will be the tightest possible, so that in no other arrangement 
more spheres could be packed into the same container''. 
\end{abstract}

\maketitle

\noindent AMS Mathematics Subject Classification 2000 (MSC2000): 52C17; 01A45, 05B40


\markboth{{\sc Achill Sch\"urmann}}{{\sc On packing spheres into containers}}

\section{Introduction}

How many equally sized spheres 
can be packed into a given container?
In 1611, {\sc Kepler} discussed this question
in his booklet 
\cite{kepler-1611} and came to the following conclusion:

\begin{quote}
``Coaptatio fiet arctissima, ut nullo praeterea ordine plures
\underline{globuli} \underline{in idem vas} compingi queant.''

\medskip

\noindent
``The (cubic or hexagonal close) packing will be the tightest possible, so that in no other arrangement more \underline{spheres} could be packed 
\underline{into the same container}.''
\end{quote}

In this note we want to show that 
Kepler's assertion is false for many containers 
(see Section \ref{sec:kepler-assert}, Corollary \ref{cor:corollary3}).
Even more general we show, roughly speaking, that the set of solutions to
the finite container problem (see below) in an Euclidean space 
of dimension $d\geq 2$ has no ``simple structure'' (see Definition \ref{def:lc}).

To make this precise,
we consider the Euclidean $d$-space $\R^d$ endowed with
inner product $\langle\cdot,\cdot\rangle$ and norm $|\cdot |$.
Let $B^d=\{\vec{x}\in\R^d : |\vec{x}|\leq 1\}$ denote the
(solid) unit sphere and $S^{d-1}=\{\vec{x}\in\R^d : |\vec{x}|= 1\}$
its boundary. Then a discrete set $X\subset \R^d$ is a {\em packing set} and 
defines a {\em sphere packing} 
$X+\frac{1}{2}B^d
=\{\vec{x}+\frac{1}{2}\vec{y} : \vec{x}\in X, \vec{y}\in B^d\}
$, 
if distinct elements $\vec{x},\vec{x}'\in X$ have distance 
$|\vec{x}-\vec{x}'|\geq 1$. The sphere packing 
is called {\em finite} if $X$ is of finite cardinality $|X|$.
Here we consider finite sphere packings contained in a convex body ({\em container}) $C$,
that is, a compact, convex subset of $\R^d$ with nonempty interior.
The {\em finite container problem} may be stated as follows.

\begin{problem}
Given $d\geq 2$, $n\in\N$ and a convex body
$C\subset\R^d$, determine
\[
\lambda(C,n)=
\min\{ \lambda >0 : \lambda C \supset X +\tfrac{1}{2}B^d \mbox{ a packing, $X\subset\R^d$ with $|X|=n$ }\}
\]
and packing sets $X$ attaining the minimum.
\end{problem}

Many specific instances of this container problem have been considered 
(see for example 
\cite{bezdek-1987}, \cite{bw-2000}, \cite{fodor-1999}, \cite{melissen-1997}, \cite{no-1999},\cite{packomania}, \cite{sccscg-2006}).
Independent of the particular choice of the container $C$,
solutions tend to densest infinite packing arrangements for growing $n$
(see Section \ref{sec:kepler-assert}, cf. \cite{cs-1995}). 
In dimension $2$ these packings are known to be arranged hexagonally.
Nevertheless, although close, solutions to the container problem 
are not hexagonally arranged for all sufficiently large $n$ 
and various convex disks $C$, as shown by the author
in \cite{schuermann-2002}, Theorem 9 
(cf. \cite{lg-1997} for corresponding computer experiments).
Here we show that a similar phenomenon is true in arbitrary 
Euclidean spaces of dimension $d\geq 2$.

We restrict ourselves to {\em smooth convex bodies} $C$ as containers.
That is, we assume the {\em support function} 
$h_C(\vec{u})=\sup \{\langle \vec{x},\vec{u} \rangle : \vec{x} \in C \}$
of $C$ is differentiable at all $\vec{u}\in\R^d\setminus\{\vec{0}\}$, or
equivalently, we require that $C$ has a unique {\em supporting hyperplane}
through each boundary point (see \cite{schneider-1993}, Chapter 1.7).

Our main result shows that families of packing sets with a 
``simple structure'' can not be solutions to the container problem 
if $C$ is smooth and $n$ sufficiently large.
This applies for example to the family of solutions to the 
{\em lattice restricted container problem}.
In it, we only consider packing sets 
which are isometric to a subset of some 
{\em lattice} (a discrete subgroup of $\R^d$).

\begin{theorem}
\label{thm:theorem1}
Let $d\geq 2$ and $C\subset\R^d$ a smooth convex body.
Then there exists an $n_0\in\N$, depending on $C$,
such that $\lambda(C,n)$ is not attained 
by any lattice packing set for $n\geq n_0$.
\end{theorem}


%
%
%

\section{Packing families of limited complexity}

The result of Theorem \ref{thm:theorem1} can be extended
to a more general class of packing sets. 

\begin{definition} \label{def:lc}
A family $\MF$ of packing sets in $\R^d$ is of
{\em limited complexity} (an {\em lc-family}), if
\begin{enumerate}
\item[(i)]
there exist isometries $\MI_X$, for each $X\in\MF$, such that
\begin{equation}
\label{eqn:lc-def}
\{
\vec{x}-\vec{y} : 
\vec{x}, \vec{y} \in \MI_X(X) \mbox{ and $X\in\MF$}
\}
\end{equation}
has only finitely many accumulation points in any bounded region.

\item[(ii)]
there exists a $\varrho>0$, such that 
for all $\vec{x}\in X$ with $X\in\MF$, 
every affine subspace spanned by some elements of 
$$\{
\vec{y} \in X :
\mbox{$|\vec{x}-\vec{y}| = 1$}
\}
$$
either contains $\vec{x}$ or its distance to $\vec{x}$
is larger than $\varrho$.
\end{enumerate}
\end{definition}

Condition (i) shows that point configurations within an arbitrarily large radius 
around a point are (up to isometries of $X$ and up to finitely many exceptions)
arbitrarily close to one out of finitely many possibilities.
Condition (ii) limits the possibilities for points at minimum distance further.
Note that the existence of a $\varrho>0$ in (ii)
follows if \eqref{eqn:lc-def} in (i) is finite within $S^{d-1}$.

An example of an lc-family 
in which isometries can be chosen so that \eqref{eqn:lc-def} is finite
in any bounded region, is the family of {\em hexagonal packing sets}. 
These are isometric copies of subsets of a {\em hexagonal lattice}, 
in which every point in the plane is at minimum distance $1$ 
to six others. For the hexagonal packing sets,
condition (ii) is satisfied for all $\varrho<\tfrac{1}{2}$. 
More general, isometric copies of subsets of a fixed lattice
give finite sets \eqref{eqn:lc-def} in any bounded region and satisfy (ii) 
for suitable small $\varrho>0$.
Similar is true for more general families of packing sets,
as for example for the {\em hexagonal close configurations} in dimension $3$ 
(see Section \ref{sec:kepler}). 

An example of an lc-family, in which the sets \eqref{eqn:lc-def} 
are not necessarily finite in any bounded region,
are the solutions to the lattice restricted container problem.
As shown at the end of Section \ref{sec:proofs},
condition (ii) in Definition \ref{def:lc} is nevertheless satisfied.
Thus we derive Theorem \ref{thm:theorem1}
from the following, more general result.

\begin{theorem} \label{thm:main}
Let $d\geq 2$, $C\subset\R^d$ a smooth convex body
and $\MF$ an lc-family of packing sets in $\R^d$.
Then there exists an $n_0\in\N$, depending on $\MF$ and $C$, such that 
$\lambda(C,n)$ is not attained 
by any packing set in $\MF$ for $n\geq n_0$.
\end{theorem}

Proofs are given in the next section.
In Section \ref{sec:extensions} we briefly mention some possible
extensions of Theorem \ref{thm:main}. 
In Section \ref{sec:kepler} we discuss consequences
for the quoted assertion of Kepler, if interpreted as a container problem 
(see Corollary \ref{cor:corollary3}).

\section{Proofs}
\label{sec:proofs}

\begin{figure}[ht]
\centerline{
\includegraphics[width=9cm]{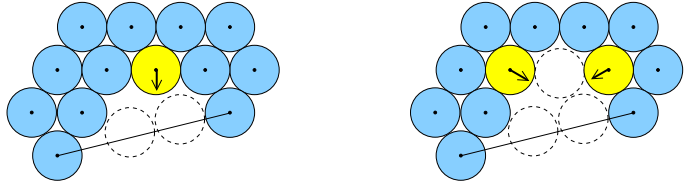}
}
\caption{Local rearrangements in a hexagonal circle packing.}
\label{fig:plane-instability}
\end{figure}

{\bf Idea. }
The proof of Theorem \ref{thm:main} is subdivided into 
four preparatory steps and corresponding propositions.
These technical ingredients are brought together at the end of this section.
Given an lc-family $\MF$ of packing sets, the idea is the following: 
We show that packing sets $X\in\MF$, 
with $|X|$ sufficiently large, allow the construction of packing sets 
$X'$ with $|X'|=|X|$ and with $X'+ \tfrac{1}{2}B^d$ 
fitting into a smaller dilate of $C$.
Roughly speaking, this is accomplished in two steps. First we show that 
``rearrangements'' of spheres near the boundary of $C$ are possible 
for sufficiently large $n$.
This allows us to obtain arbitrarily large regions 
in which spheres have no contact, respectively in which points 
of $X'$ have distance greater than $1$ to all other points  
(Proposition \ref{prop:prop2}, depending on property (i) of Definition \ref{def:lc}).
Such an initial modification then allows rearrangements
of all spheres (Proposition \ref{prop:prop3} and \ref{prop:prop4}, 
depending on property (ii) of Definition \ref{def:lc}), 
so that the resulting packing fits into a smaller dilate of $C$.
For example, consider a hexagonal packing in the plane:
It is sufficient to initially rearrange (or remove) two disks in order
to subsequently rearrange all other disks, so that no disk 
is in contact with others afterwards 
(see Figure \ref{fig:plane-instability}, cf. \cite{schuermann-2002}).

How do we know that the new sphere packings $X'+\tfrac{1}{2}B^d$ fit
into a smaller dilate of $C$? Consider 
$$
\lambda(C,X)
=
\min
\{
\lambda > 0 :
\lambda C\supset \vec{t}+X +\tfrac{1}{2}B^d
\mbox{ for some } \vec{t}\in\R^d
\}
$$
for a fixed finite packing set $X$. 
Here and in the sequel we use
$\vec{t}+X$ to abbreviate $\{\vec{t}\}+X$.
Clearly 
$$\lambda(C,n)=\min\{ \lambda(C,X) : \mbox{ $X$ is a packing set with $|X|=n$ }\},$$
and $\lambda(C,X')<\lambda(C,X)$
whenever the {\em convex hull} $\conv X'$ of $X'$
(and hence $X'$ itself) is contained in the interior
$\interior\conv X$ of the convex hull of $X$.
Thus in order to prove that $X$ does not attain $\lambda(C,|X|)$
for any convex container $C$, it is sufficient to describe a way
of attaining a packing set $X'$ with $|X'|=|X|$ and
\begin{equation} \label{eq:suf-cond}
X' \subset \interior\conv X
.
\end{equation}

\medskip

{\bf I. } 
Let us first consider the ``shapes''
of packing sets $X_n$ attaining $\lambda(C,n)$.
Here and in what follows, $X_n$ denotes
a packing set with $|X_n|=n$.

In order to define the ``shape'', let 
\[
R(M)=\min\{ R \geq 0 : M\subset \vec{t} + R B^d \mbox{ for some }\vec{t}\in\R^d\}
\]
denote the circumradius of a compact set $M\subset\R^d$ and 
let $\vec{c}(M)$ denote the center of its circumsphere.
Hence $M\subseteq \vec{c}(M)+R(M)B^d$.
Then the {\em shape} of $M$ is defined by 
\[
\MS(M)=\left(\conv(M)-\vec{c}(M)\right) / R(M) \subset B^d
.
\]
The family of nonempty compact subsets in $\R^d$ can be
turned into a metric space, for example with the {\em Hausdorff metric} 
(cf. \cite{schneider-1993}).
Shapes of packing sets $X_n$ attaining $\lambda(C,n)$ converge to the shape of $C$,
that is,
\begin{equation}    \label{shape-convergence}
\lim_{n\to\infty} \MS\left(X_n\right) = \MS(C)
.
\end{equation}
This is seen by ``reorganizing elements'' in a hypothetical convergent 
subsequence of $\{X_n\}_{n\in\N}$ not satisfying \eqref{shape-convergence}.

The convergence of shapes 
leads for growing $n$ to shrinking sets of
{\em outer (unit) normals} 
\begin{equation} \label{eqn:normals}
\{ \vec{v}\in S^{d-1} : \langle \vec{v}, \vec{x} \rangle \geq \langle \vec{v}, \vec{y} \rangle 
 \mbox{ for all } \vec{y}\in \conv X_n\}
\end{equation}
at boundary points $\vec{x}$ of the {\em center polytope} $\conv X_n$.
For general terminology and results on convex polytopes 
used here and in the sequel we refer to \cite{ziegler-1998}.

Since $C$ is smooth, the sets of outer normals \eqref{eqn:normals} at boundary points
of $\conv X_n$ become uniformly small for large $n$. 
Also, within a fixed radius around a boundary point, 
the boundary of $\conv X_n$ becomes ``nearly flat'' for growing $n$.

\begin{proposition} \label{prop:prop1}
Let $d\geq 2$ and $C\subset \R^d$ a smooth convex body.
Let $\{X_n\}$ be a sequence of packing sets in $\R^d$ attaining $\lambda(C,n)$. 
Then 
\begin{enumerate}
\item[(i)]
for $\varepsilon > 0$ 
there exists an $n_1\in\N$, depending on $C$ and $\varepsilon$, 
such that for all $n\geq n_1$, outer normals $\vec{v},\vec{v}'\in S^{d-1}$ of
$\conv X_n$ at $\vec{x}\in X_n$ satisfy
$$
|\vec{v}-\vec{v}'|
<\varepsilon
;
$$ 

\item[(ii)]
for $\varepsilon>0$ and $r>0$ there exists an $n_1\in\N$, depending on $C$, $\varepsilon$ and $r$, 
such that for all $n\geq n_1$, and for $\vec{x},\vec{x}'\in \bd \conv X_n$ with $|\vec{x}-\vec{x}'|\leq r$,
outer normals $\vec{v}\in S^{d-1}$ of $\conv X_n$ at $\vec{x}$ satisfy
$$
\langle \vec{v}, \vec{x}-\vec{x}' \rangle > -\varepsilon
.
$$
\end{enumerate}
\end{proposition}

\medskip

{\bf II. }
In what follows we use some additional terminology.
Given a packing set $X$, we say 
$\vec{x}\in X$ is in a {\em free position},
if the set 
\[
\MN_X(\vec{x})=\{ \vec{y}\in X : |\vec{x}-\vec{y}|=1 \}
\]
is empty.
If some $\vec{x}\in X$ is not 
contained in 
$\interior \conv \MN_X(\vec{x})$, then it is possible to 
obtain a packing set $X'=X\setminus\{\vec{x}\}\cup \{\vec{x}'\}$
in which $\vec{x}'$ is in a free position. 
We say $\vec{x}$ is {\em moved to a free position} in this case (allowing $\vec{x}'=\vec{x}$). 
We say $\vec{x}$ is {\em moved into or within a set} $M$ (to a free position), 
if $\vec{x}'\in M$. 
Note, in the resulting packing set $X'$
less elements may have minimum distance $1$ to others,
and therefore possibly further elements can be moved to free positions.

Assuming $X\in\MF$ attains $\lambda(C,|X|)$ with $|X|$ sufficiently large,
the following proposition shows that it is possible to move elements of $X$ 
into free positions within an arbitrarily large region,
without changing the center polytope $\conv X$.

\begin{proposition} \label{prop:prop2}
Let $d\geq 2$ and $R>0$. Let $C\subset \R^d$ a smooth convex body and
$\MF$ a family of packing sets in $\R^d$ satisfying (i) of Definition \ref{def:lc}.
Then there exists an $n_2\in\N$, depending on $R$, $\MF$ and $C$, 
such that for all $X \in\MF$ attaining $\lambda(C,|X|)$
with $|X|\geq n_2$, there exists a $\vec{t}_{X}\in\R^d$ with 
\begin{enumerate}
\item[(i)]
$(\vec{t}_{X}+RB^d)\subset \conv X$, and 
\item[(ii)]
all elements of $X\cap\interior (\vec{t}_{X}+RB^d)$ can be moved to free positions
by subsequently moving elements of $X\cap \interior \conv X$ 
to free positions within $\interior \conv X$.
\end{enumerate}
\end{proposition}

\begin{proof}
{\em Preparations. }
By applying suitable isometries to the packing sets in $\MF$
we may assume that 
\begin{equation}  \label{eqn:simpl-asump}
\{
\vec{y} : 
\mbox{ $\vec{y} \in X - \vec{x}$ with $|\vec{y}| < r$ for $\vec{x}\in X$ and $X\in\MF$}
\}
\end{equation}
has only finitely many accumulation points for every $r>1$. 
For each $X$, the container $C$ is transformed to possibly 
different isometric copies.
This is not a problem though, since the container 
is not used aside of Proposition \ref{prop:prop1},
which is independent of the chosen isometries.
Note that the smoothness of $C$ is implicitly used here.

We say $\vec{x}\in X$ is {\em moved in direction $\vec{v}\in S^{d-1}$},
if it is replaced by an $\vec{x}'$ on the ray 
$\{\vec{x}+\lambda \vec{v} : \lambda \in\R_{> 0} \}$.
Note that it is possible to move $\vec{x}$ in direction $\vec{v}\in S^{d-1}$ 
to a free position, if 
\begin{equation}  \label{eqn:lattice-blocking-sets}
\MN_X(\vec{x},\vec{v}) 
= \{ \vec{w} \in \MN_{X}(\vec{x})-\vec{x} : \langle \vec{v}, \vec{w} \rangle > 0 \}
\end{equation}
is empty.
If we want a fixed $\vec{x}\in X$ to be moved to a free position,
in direction $\vec{v}\in S^{d-1}$ say, we have to move the elements 
$\vec{y}\in \vec{x}+\MN_X(\vec{x},\vec{v})$ first. 
In order to do so,
we move the elements of $\vec{y}+\MN_X(\vec{y},\vec{v})$ to free positions, and so on.
By this we are lead to the definition of the {\em access cone} 
\begin{equation}  \label{eqn:access_cones}
\mbox{acc}_{\MF,n}(\vec{v})
=
\mbox{pos}\left\{
\MN_{X}(\vec{x},\vec{v})  :  \vec{x}\in X \mbox{ for } X \in \MF \mbox{ with } |X|\geq n 
\right\}
\end{equation}
of $\MF$ and $n$ in direction $\vec{v}\in S^{d-1}$. 
Here, 
$$
\mbox{pos}(M) = \{ \sum_{i=1}^m \lambda_i\vec{x}_i : m\in\N, \lambda_i\geq 0 \mbox{ and } \vec{x}_i\in M \mbox{ for $i=1,\dots,m$ } \}
$$ 
denotes the {\em positive hull} of a set $M\subset \R^d$, which is by definition a convex cone.  
Note that $\mbox{acc}_{\MF,n}(\vec{v})$ is contained in
the halfspace $\{\vec{x} \in\R^d : \langle \vec{v}, \vec{x} \rangle \geq 0 \}$
and that $\mbox{acc}_{\MF,n}(\vec{v}) \subseteq \mbox{acc}_{\MF,n'}(\vec{v})$
whenever $n\geq n'$.

By the assumption that
\eqref{eqn:simpl-asump} has only finitely many accumulation points for $r>1$,
there exist only finitely many limits 
$\lim_{n\to\infty}\left(\mbox{acc}_{\MF,n}(\vec{v})\cap B^d\right)$. 
Here, limits are defined using the Hausdorff metric on the set 
of nonempty compact subsets of $\R^d$ again.

\smallskip

{\em Strategy. }
We choose a $\vec{v}\in S^{d-1}$ such that there exists an $\varepsilon>0$
with 
$$
\lim_{n\to\infty}\left(\mbox{acc}_{\MF,n}(\vec{v})\cap B^d\right)
=
\lim_{n\to\infty}\left(\mbox{acc}_{\MF,n}(\vec{v}')\cap B^d\right)
,
$$ 
for all $\vec{v}'$ in the $\varepsilon$-{\em neighborhood}
$S_{\varepsilon}(\vec{v})=S^{d-1}\cap (\vec{v}+{\varepsilon} B^d)$
of $\vec{v}\in S^{d-1}$.

In order to prove the proposition, 
we show the following for every $X\in\MF$, attaining $\lambda(C,|X|)$
with $|X|$ sufficiently large:
There exists a $\vec{t}_{X}\in\R^d$ such that
\begin{enumerate}
\item[(i')]
$(\vec{t}_{X}+RB^d) + \mbox{acc}_{\MF,n}(\vec{v})$ does not intersect $X \cap \bd \conv X$, while

\item[(ii')]
$(\vec{t}_{X}+RB^d) \subset  \conv X$.
\end{enumerate}
It follows that $\bd \conv X$ has to intersect 
the unbounded set 
\begin{equation} \label{eqn:acc-translate}
(\vec{t}_{X}+RB^d) + \mbox{acc}_{\MF,n}(\vec{v})
\end{equation}
and by the definition of the access cone it is possible
to move the elements in $X\cap \interior(\vec{t}_X+RB^d)$ to free positions as asserted.
For example,
after choosing a direction $\vec{v}'\in S_{\varepsilon}(\vec{v})$,
we may subsequently pick 
non-free elements $\vec{x}$ in \eqref{eqn:acc-translate} 
with maximal $\langle \vec{x},\vec{v}' \rangle$.
These elements can be moved to a free position within $\interior \conv X$,
since $\MN_X(\vec{x},\vec{v}')$ is empty by the definition
of the access cone. 

\smallskip

{\em Bounding the boundary intersection. } 
We first estimate the size of the intersection
of \eqref{eqn:acc-translate} with $\bd \conv X$.  
For $\vec{v}'\in S_{\varepsilon}(\vec{v})$ and $n\in\N$, we consider the sets 
$$
M(\vec{v}',n)=\{ \vec{x}\in RB^d+\mbox{acc}_{\MF,n}(\vec{v}) : \langle \vec{x} , \vec{v}' \rangle = R\}
.
$$
By the definition of the access cones \eqref{eqn:access_cones},  
$M(\vec{v}',n)\subseteq M(\vec{v}',n')$ for $n\geq n'$.
We choose 
$$
r
>
\sup
\{
|\vec{x}-\vec{y} | : \vec{x},\vec{y}\in M(\vec{v}',n) 
\mbox{ with } \vec{v}'\in S_{\varepsilon}(\vec{v})
\}
,
$$
as a common upper bound on the 
diameter of the sets $M(\vec{v}',n)$ with $n$ sufficiently large, 
say $n\geq n'$. 
Note that $R$ as well as $\MF$, $\vec{v}$ and $\varepsilon$  
have an influence on the size of $r$ and $n'$.

By Proposition \ref{prop:prop1} (ii)
we can choose $n'$ possibly larger 
to ensure the following for all $X\in\MF$ attaining $\lambda(C,|X|)$
with $|X|\geq n'$:
{\em 
The intersection of \eqref{eqn:acc-translate} with $\bd\conv X$
has a diameter less than $r$, no matter which 
$\vec{t}_{X}\in \conv X$ 
at distance $R$ to $\bd\conv X$ we choose. 
Moreover, $(\vec{t}_{X}+RB^d) \subset \conv X$.
}

\smallskip

{\em Ensuring an empty intersection. } 
It remains to show that for $X\in\MF$, attaining $\lambda(C,|X|)$ 
with $|X|$ sufficiently large,
$\vec{t}_{X}$ can be chosen such that \eqref{eqn:acc-translate} does not
intersect $X\cap \bd \conv X$.
For this we prove the following
{claim:} 
{\em 
There exists an $n''$, depending on $r$, $\vec{v}$ and $\varepsilon$, such that
for all $X\in\MF$ with $|X|\geq n''$, there exists
a vertex $\vec{x}$ of $\conv X$ with outer normal $\vec{v}'\in S_{\varepsilon}(\vec{v})$
and 
\begin{equation} \label{eqn:final-cond}
\{\vec{x}\}= X \cap (\bd \conv X ) \cap (\vec{x}+r B^d)
.
\end{equation}
}

Thus these vertices $\vec{x}$ have a distance larger than $r$ 
to any other element of $X \cap \bd \conv X$.
Therefore, by choosing $n_2\geq \max \{n', n''\}$,
we can ensure that there exists a $\vec{t}_{X}\in\R^d$
at distance $R$ to $\bd \conv X$ 
such that (i') and (ii') are satisfied for all $X\in\MF$ 
attaining $\lambda(C,|X|)$ with $|X| \geq n_2$.
Note that $n'$, $n''$, and hence $n_2$, 
depend on the choice of $\vec{v}$ and $\varepsilon$.
But we may choose $\vec{v}$ and $\varepsilon$,
depending on $\MF$, so that $n_2$ can be chosen as small as possible.
In this way we get an $n_2$ which solely depends on $R$, $\MF$ and $C$.

It remains to prove the claim.
Since \eqref{eqn:simpl-asump} has only finitely many accumulation points, 
the set of normals $\vec{v}' \in S^{d-1}$ with hyperplane 
$\{ \vec{y} \in\R^d : \langle \vec{v}',\vec{y}\rangle = 0\}$
running through $\vec{0}$ and an accumulation point 
$\vec{y}$ of \eqref{eqn:simpl-asump}
all lie in the union $\MU_r$ of finitely many linear subspaces of dimension $d-1$.
Thus for any $\delta>0$ the normals of these hyperplanes 
all lie in $\MU_{r,\delta} = \MU_r+\delta B^d$ if we choose $|X|$ 
sufficiently large, depending on $\delta$. 
By choosing $\delta$ small enough, we 
find a $\vec{v}'\in S_{\varepsilon}(\vec{v})$ 
with $\vec{v}'\not\in \MU_{r,\delta}$.
Moreover, there exists an $\varepsilon'>0$ such that
$S_{\varepsilon'}(\vec{v}')\cap \MU_{r,\delta} = \emptyset$.
Since every center polytope $\conv X$ has a vertex $\vec{x}$ 
with outer normal $\vec{v}'$,
we may choose $|X|$ sufficiently large by Proposition \ref{prop:prop1} (i)
(applied to $2 \varepsilon '$),
such that $\conv X$ has no outer normal in $\MU_{r,\delta}$ at $\vec{x}$.

Moreover, for sufficiently large $|X|$, faces of $\conv X$ intersecting
$\vec{x}+rB^d$ can not contain any vertex in $X\cap (\vec{x}+rB^d)$ 
aside of $\vec{x}$. Thus by construction, there exists an $n''$
such that \eqref{eqn:final-cond} holds for all $X\in\MF$ with $|X|\geq n''$.
This proves the claim and therefore the proposition.
\end{proof}

Note that the proof offers the possibility 
to loosen the requirement on $\MF$ a bit,
for the prices of introducing another parameter: 
For suitable large $r$, depending on $\MF$,
the proposition holds, if 
\begin{enumerate}
\item[(i')]
there exist isometries $\MI_X$ for each $X\in\MF$, such that
$$
\{
\vec{x}-\vec{y} : 
\vec{x}, \vec{y} \in \MI_X(X) \mbox{ and $X\in\MF$ }
\}
$$
has only finitely many accumulation points within $r B^d$.
\end{enumerate}

\medskip

{\bf III. }
For all $X\in\MF$ attaining $\lambda(C,|X|)$, with $|X|$ sufficiently large,
we are able to obtain {\em contact free regions} 
$(\vec{t}_{X}+RB^d)\subset \conv X$, with $R$ as large as we want,
by Proposition \ref{prop:prop2}.
That is, we can modify these packing sets $X$ by moving elements to free positions
within $\interior (\vec{t}_X+RB^d)$.
By choosing $R$ large enough,
such an initial contact free region allows to move further elements 
to free positions. The following proposition takes care of interior points.

\begin{proposition} \label{prop:prop3}
Let $d\geq 2$ and $\MF$ a family of packing sets in $\R^d$ 
satisfying (ii) in Definition \ref{def:lc} with $\varrho>0$. 
Let $R\geq \tfrac{1}{\varrho}$, $X\in\MF$ and $\vec{x}\in X \cap \interior \conv X$.
Let $\vec{t}\in\R^d$ with $|\vec{t}-\vec{x}|\leq R+\tfrac{\varrho}{2}$
and with all elements of $X \cap (\vec{t}+RB^d)$ in a free position.
Then $\vec{x}$ can be moved to a free position within $\interior \conv X$.
\end{proposition}

\begin{proof}
Assume $\vec{x}\in\interior \conv \MN_X(\vec{x})$.
By the assumption on $\MF$,  
$$
\vec{x}+\varrho B^d\subset\interior \conv \MN_X(\vec{x})
.
$$
Thus there exists a $\vec{y}\in \MN_X(\vec{x})$, such that the
orthogonal projection $\vec{y}'$ of $\vec{y}$ onto the line through
$\vec{x}$ and $\vec{t}$ satisfies $|\vec{y}'-\vec{x}|\geq \varrho$
and $|\vec{y}'-\vec{t}|\leq R-\tfrac{\varrho}{2}$.
Then 
$$
|\vec{y}-\vec{t}|^2 
= 
|\vec{y}'-\vec{t}|^2 + |\vec{y}-\vec{y}'|^2
\leq \left(R-\tfrac{\varrho}{2}\right)^2 + \left(1 - \varrho^2\right) < R^2
.
$$
Thus $\vec{y}$ is in a free position by the assumptions of the 
proposition, which contradicts $\vec{y}\in \MN_X(\vec{x})$.
\end{proof}

\medskip

{\bf IV. }
After Propositions \ref{prop:prop2} and \ref{prop:prop3}
it remains to take care of points in $X \cap \bd \conv X$,
for $X\in\MF$ attaining $\lambda(C,|X|)$, and with $|X|$ sufficiently large.
It turns out that these points can all be moved to free positions
within $\interior \conv X$. As a consequence we obtain the following.

\begin{proposition} \label{prop:prop4}
Let $d\geq 2$, $C\subset \R^d$ a smooth convex body
and $\MF$ a family of packing sets in $\R^d$ satisfying (ii) of Definition \ref{def:lc}.
Then there exists an $n_4\in\N$, depending on $C$ and $\MF$,
such that $X\in\MF$ with $|X|\geq n_4$ does not attain $\lambda(C,|X|)$,
if all elements of $X\cap \interior \conv X$ are in a free position.
\end{proposition}

\begin{proof}
Let $\varrho>0$ as in (ii) of Definition \ref{def:lc}.
We choose $n_4$ by Proposition \ref{prop:prop1} (ii), applied to $\varepsilon =\varrho$ and $r=1$.
Assume $X\in\MF$ with $|X|\geq n_4$ attains $\lambda(C,|X|)$
and all elements of $X\cap \interior \conv X$ are in a free position.
We show that every element $\vec{x}\in X\cap \bd \conv X$
can be moved to a free position into $\interior \conv X$.
This gives the desired contradiction,
because after moving (in an arbitrary order) all $X \cap \bd \conv X $ 
to free positions into $\interior \conv X$, we
obtain a packing set $X'$ with $|X'|=|X|$
and $X'\subset \interior \conv X$. 

It is possible to move a given $\vec{x}\in X\cap \bd \conv X$ 
to a free position $\vec{x}'=\vec{x}+\delta \vec{v}$
for a (sufficiently small) $\delta>0$, if $\vec{v}\in S^{d-1}$ is contained in 
the non-empty polyhedral cone
$$
C_{\vec{x}}
=
\left\{
\vec{v}\in \R^d:
\langle \vec{v}, \vec{y}-\vec{x}\rangle \leq 0 
\mbox{ for all } \vec{y}\in\MN_X(\vec{x})
\right\}
.
$$
If $\vec{v}\in C_{\vec{x}}$ can be chosen, so that $\vec{x}'\in\interior \conv X$,
the assertion follows.
Otherwise, because $C_{\vec{x}}$ and $\conv X$ are convex, there exists a hyperplane 
through $\vec{x}$, with normal $\vec{w}\in S^{d-1}$, 
which separates $\conv X$ and $\vec{x}+C_{\vec{x}}$.
That is, 
we may assume that 
$$
\vec{w}
\in 
\pos \left\{
\vec{y}-\vec{x} : \vec{y} \in\MN_X(\vec{x})
\right\}
$$
and $-\vec{w}$ is an outer normal of $\conv X$ at $\vec{x}$.

Then for some $\delta>0$, there exists a point $\vec{z} = \vec{x} + \delta \vec{w} \in \bd \conv \MN_X(\vec{x})$,
which is a convex combination of some $\vec{y}_1,\dots,\vec{y}_k \in \MN_X(\vec{x})$.
That is, there exist $\alpha_i\geq 0$ with 
$
\sum_{i=1}^k \alpha_i = 1 
$
and
$
\vec{z}=\sum_{i=1}^k \alpha_i \vec{y}_i
$
.
Therefore 
$$
\delta = \langle \vec{z}-\vec{x} , \vec{w} \rangle
= \sum_{i=1}^k \alpha_i \langle \vec{y}_i -\vec{x} , \vec{w} \rangle
< \varrho,
$$
because $\langle \vec{y}_i - \vec{x}, \vec{w} \rangle < \varrho$
due to $|X|\geq n_4$ and $\vec{y}_i\in\bd \conv X$.
This contradicts the assumption on $\MF$ with respect to $\varrho$ though.
\end{proof}

\medskip

{\bf Finish. } 
The proof of Theorem \ref{thm:main}
reduces to the application of Propositions \ref{prop:prop1}, \ref{prop:prop2},
\ref{prop:prop3} and \ref{prop:prop4}.
Let $\MF$ be an lc-family of packing sets in $\R^d$,
with a $\varrho>0$ as in (ii) of Definition \ref{def:lc}.
We choose $R \geq 1/\varrho$ and $n_2$ and $n_4$ 
according to Propositions \ref{prop:prop2} and \ref{prop:prop4}.
By Proposition \ref{prop:prop1} (ii), we choose $n_1$ such that packing sets 
$X$ attaining $\lambda(C,|X|)$ with $|X|\geq n_1$ satisfy the following:
For each $\vec{x}\in X$, there exists a $\vec{t}\in \R^d$ with 
$|\vec{x}-\vec{t}| = R + \tfrac{\varrho}{2}$ and 
$\vec{t}+ R B^d \subset \conv X$.

We choose $n_0\geq \max\{n_1,n_2,n_4\}$
and assume that $X\in\MF$ with $|X|\geq n_0$ attains $\lambda(C,|X|)$. 
By Proposition \ref{prop:prop2} 
we can modify the packing set $X$ to obtain a new packing set $X'$
with a contact free region $(\vec{t}_{X}+R B^d)\subset \interior \conv X$,
and with the same points 
$X'\cap \bd \conv X'=X\cap \bd \conv X$ on the boundary of 
the center polytope $\conv X' = \conv X$.

The following gives a possible order, in which 
we may subsequently move non-free elements $\vec{x}\in X \cap \interior\conv X$
to free positions:
By the choice of $n_0$ we can guarantee that for each $\vec{x}\in X \cap \interior\conv X$,
there exists a $\vec{t}$ with $|\vec{x}-\vec{t}|\leq R+\tfrac{\varrho}{2}$
and $\vec{t}+ R B^d \subset \conv X$. Let $\vec{t}_{\vec{x}}$ be the
$\vec{t}$ at minimal distance to $\vec{t}_X$.
Then among the non-free $\vec{x}\in \interior\conv X$, the one with 
minimal distance $|\vec{t}_{\vec{x}}-\vec{t}_X|$ satisfies 
the assumptions of Proposition \ref{prop:prop3},
because a non-free element $\vec{y}\in X \cap (\vec{t}_{\vec{x}}+B^d)$
would satisfy $|\vec{t}_{\vec{y}}-\vec{t}_X|<|\vec{t}_{\vec{x}}-\vec{t}_X|$ 
due to $\conv\{\vec{t}_{\vec{x}},\vec{t}_X\}+B^d\subset \conv X$.

Thus by Proposition \ref{prop:prop3} we can subsequently move the non-free elements 
within $X\cap \interior\conv X$ to free positions. 
By this we obtain a contradiction to Proposition \ref{prop:prop4}, 
which proves the theorem.

\medskip

{\bf The lattice packing case. }
We end this section with the proof of Theorem \ref{thm:theorem1}. 
We may apply Theorem \ref{thm:main} after showing 
that the family of solutions to the lattice restricted container problem 
is of limited complexity.
The space of lattices can be turned into a topological space
(see \cite{gl-1987}).  
The convergence of a sequence $\{\Lambda_n\}$ of lattices 
to a lattice $\Lambda$ in particular involves that sets of lattice points 
within radius $r$ around a lattice point tend to 
translates of $\Lambda\cap rB^d$ for growing $n$. 
As a consequence, a convergent sequence of packing lattices,
as well as subsets of them, form an lc-family.
Solutions to the lattice restricted container problem tend for growing $n$ towards subsets 
of translates of {\em densest packing lattices} (see \cite{zong-1999}). 
These lattices are the solutions of the {\em lattice (sphere) packing problem}.
Up to isometries, there exist only finitely 
many of these lattices in each dimension (see \cite{zong-1999}). 
Thus the assertion follows, since a finite union of lc-families is
an lc-family.

\section{Extensions}

\label{sec:extensions}

Let us briefly mention some possible extensions of Theorem \ref{thm:main}.
These have been treated in \cite{schuermann-2002}
for the $2$-dimensional case and could be directions for further 
research.

\medskip

{\em Packings of other convex bodies. }  
Instead of sphere packings, we may consider
packings $X+K$ for other convex bodies $K$.
If the {\em difference body} $DK=K-K$ is strictly convex,
then the proofs can be applied after some modifications:
Instead of measuring distances with the norm $|\cdot|$ given
by $B^d$, we use the
norm $|\vec{x}|_{DK}=\min \{\lambda>0 : \lambda \vec{x} \in DK\}$
given by $DK$. The strict convexity of $DK$ is then used for the key fact, 
that elements $\vec{x}$ of a packing set $X$ can be moved to a free
position, whenever they are not contained in $\interior \conv \MN_X(\vec{x})$
(see {\bf II} in Section \ref{sec:proofs}). 
Note though that the sets in  
\eqref{eqn:lattice-blocking-sets} and depending definitions have 
to be adapted for general convex bodies.

\medskip

{\em Packings in other containers. }
The restriction to smooth convex containers simplifies the proof,
but we strongly believe that Theorem \ref{thm:main} is valid for other
containers as well, e.g. certain polytopes.
On the other hand there might exist containers for which
Theorem \ref{thm:main} is not true.
In particular in dimension $3$ it seems very likely that Theorem \ref{thm:main}
is not true for polytopal containers $C$ with all their facets lying 
in planes containing hexagonal sublattices of the {\em fcc lattice} (see Section \ref{sec:kepler-assert}). 
That is, for these polytopal containers $C$ we conjecture the existence of infinitely many $n$,
for which subsets of the {\em fcc lattice} attain $\lambda(C,n)$. 
An example for at least ``local optimality'' 
of sphere packings (with respect to differential perturbations) 
in suitable sized tetrahedra was given by 
Dauenhauer and Zassenhaus \cite{dz-1987}.
A proof of ``global optimality'' seems extremely difficult though, as it would provide a new
proof of the sphere packing problem (``Kepler conjecture'', see Section \ref{sec:kepler}).

\medskip

{\em Other finite packing problems. }
Similar ``phenomena'' occur for other packing problems. 
For example, if we consider finite packing sets $X$ with 
minimum diameter or surface area of $\conv X$, or maximum 
parametric density with large parameter
(cf. \cite{cfg-1991}, \cite{bhw-1994}, \cite{boeroeczky-2004}, \cite{bp-2005}).
This is due to the fact that the shapes of solutions tend to
certain convex bodies, e.g. a sphere.

\section{Kepler's assertion} \label{sec:kepler-assert}

\label{sec:kepler}

Kepler's statement, quoted in the introduction, 
was later referred to as the origin of the 
famous sphere packing problem known as the {\em Kepler conjecture}
(cf. e.g. \cite{hales-2002} p.5, \cite{hsiang-2001} p.4).
In contrast to the original statement, this problem 
asks for the maximum sphere packing density (see \eqref{eq:delta} below) of an infinite arrangement
of spheres, where the ``container'' is the whole Euclidean space.
As a part of Hilbert's famous problems \cite{hilbert-1901},
it attracted many researchers in the past. 
Its proof by Hales with contributions of Ferguson 
(see \cite{hales-2002}, \cite{hales-2005}, \cite{hales-2006}), 
although widely accepted, had been a matter of discussion 
(cf. \cite{lagarias-2002}, \cite{szpiro-2003}, \cite{fl-2006}).

\vspace*{0.1cm}
\noindent
\begin{minipage}{8.5cm}
\hspace*{0.3cm}
Following Kepler \cite{kepler-1611}, 
the {\em cubic or hexagonal close packings} in $\R^3$ can
be described via two dimensional layers of spheres, 
in which every sphere center belongs to a planar square grid,
say with minimum distance $1$. 
These layers are stacked (in a unique way) such that 
each sphere in a layer touches exactly four spheres of the layer above
and four of the layer below.  
\end{minipage}
\hfill
\begin{minipage}{4cm}
\hfill
\includegraphics[width=3.5cm]{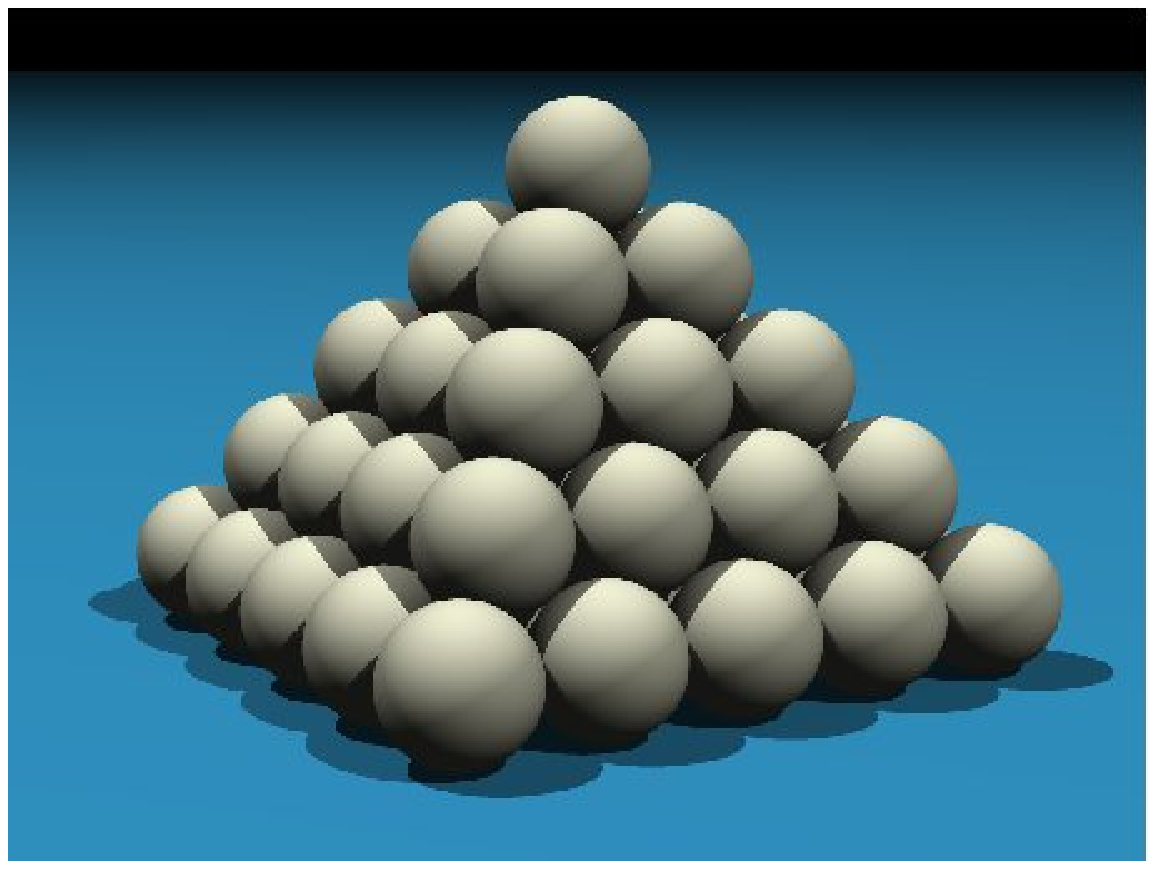}
\end{minipage}
\vspace*{0.03cm}

The packing attained in this way is the well known
{\em face centered cubic (fcc) lattice packing}. 
We can build up the fcc lattice by planar hexagonal layers as well,
but then there are two choices for each new layer 
to be placed, and only one of them yields an fcc lattice packing. 
All of them, including the uncountably many non-lattice packings, 
are referred to as {\em hexagonal close packings (hc-packings)}.
Note that the family of hc-packings 
is of limited complexity, because up to isometries
they can be built from a fixed hexagonal layer.

Let
\[
n(C)=\max\{ |X| : C \supset X +\tfrac{1}{2}B^d
            \mbox{ is a packing } \}.
\]
Then in our terminology Kepler asserts that,  
in $\R^3$, $n(C)$ is attained by hc-packings.
His assertion, if true, would imply an ``answer'' to the sphere packing problem
(Kepler conjecture), 
namely that the density of the densest infinite sphere packing
\begin{equation} \label{eq:delta}
\delta_d=
\displaystyle \limsup_{\lambda \to \infty} 
\frac{n(\lambda C) \cdot \vol(\frac{1}{2}B^d)}{\vol(\lambda C)}
\end{equation}
is attained by hc-packings for $d=3$; hence 
$\delta_3=\pi/\sqrt{18}$.
Note that this definition of density is independent
of the chosen convex container $C$ (see \cite{hlawka-1949} or \cite{gl-1987}).

As a consequence of Theorem \ref{thm:main}, Kepler's assertion turns out to be false,
even if we think of arbitrarily large containers.
Consider for example the containers $\lambda(C,n)C$ for $n\geq n_0$.

\begin{corollary}
\label{cor:corollary2}
Let $d\geq 2$, $C\subset\R^d$ a smooth convex body
and $\MF$ an lc-family of packing sets in $\R^d$.
Then there exist arbitrarily 
large $\lambda$ such that $n(\lambda C)$ 
is not attained by packing sets in $\MF$. 
\end{corollary}

We may as well think of arbitrarily small spheres
packed into a fixed container $C$.
For $r>0$, we call $X + r B^d$ a sphere packing if distinct
elements $\vec{x}$ and $\vec{x}'$ of $X$ have distance
$|\vec{x}-\vec{x}'|\geq 2r$. 
Specializing to $\R^3$, the following corollary of
Theorem \ref{thm:main} refers directly to Kepler's assertion.

\begin{corollary}
\label{cor:corollary3}
Let $C\subset\R^3$ a smooth convex body.
Then there exist arbitrarily small $r>0$, such that
$$
\max\{ |X| : C \supset X + r B^d \mbox{ is a packing } \}
$$
is not attained by fcc or hexagonal close packing sets.
\end{corollary}

\section*{Acknowledgments}

I like to thank Thomas C. Hales, Tyrrell B. McAllister, Frank Vallentin, 
J\"org M. Wills, G\"unter M. Ziegler and the two anonymous referees
for many helpful suggestions.

\providecommand{\href}[2]{#2}


\begin{thebibliography}{UBW94}

\bibitem[BHW94]{bhw-1994}
  U.~Betke, M.~Henk and J.M. Wills, \emph{Finite and infinite packings}, J. Reine Angew. Math. \textbf{453} (1994), 165--191.

\bibitem[Bez87]{bezdek-1987}
K.~Bezdek, \emph{Densest packing of small number of congruent spheres in polyhedra},
  Ann. Univ. Sci. Budapest E\"otv\"os Sect. Math., \textbf{30} (1987), 177--194.

\bibitem[B\"or04]{boeroeczky-2004}
K.~B\"or\"oczky Jr., \emph{Finite packing and covering}, Cambridge university
  press, Cambridge, 2004.

\bibitem[BW04]{bw-2000}
K.~B\"or\"oczky Jr., G.~Wintsche, \emph{Sphere packings in the regular crosspolytope},
  Ann. Univ. Sci. Budapest E\"otv\"os Sect. Math., \textbf{43} (2000), 151--157.

\bibitem[BP05]{bp-2005}
P.~Brass and J.~Pach, \emph{Problems in discrete geometry}, Springer, New York, 2005.

\bibitem[CS95]{cs-1995}
  J.H.~Conway and N.J.A.~Sloane,
  \emph{What are all the best sphere packings in low dimensions?},
  Discrete Comput. Geom., \textbf{13} (1995), 282--403.

\bibitem[DZ87]{dz-1987}
M.H. Dauenhauer and H.J. Zassenhaus, \emph{Local optimality of the critical
  lattice sphere-packing of regular tetrahedra}, Discrete Math. \textbf{64}
  (1987), 129--146.

\bibitem[FCG91]{cfg-1991}
K.J.~Falconer H.T.~Croft and R.K. Guy, \emph{Unsolved problems in geometry},
  Springer, New York, 1991.

\bibitem[FL06]{fl-2006}
G.~Fejes T\'oth and J.C.~Lagarias,
  \emph{Guest Editors Foreword},
  Discrete Comput. Geom., \textbf{36} (2006), 1--3.

\bibitem[Fod99]{fodor-1999}
F.~Fodor, \emph{The densest packing of 19 congruent circles in a circle},
  Geom. Dedicata, \textbf{74} (1999), 139--145.

\bibitem[GL87]{gl-1987}
P.M.~Gruber and C.G. Lekkerkerker, \emph{Geometry of numbers}, North-Holland,
  Amsterdam, 1987.

\bibitem[Hal02]{hales-2002}
T.C.~Hales, \emph{An overview of the {K}epler conjecture}, 16 pages, first
  in a series, \href{http://www.arxiv.org/math/9811071}{arXiv:math.MG/9811071} (v2).

\bibitem[Hal05]{hales-2005}
T.C.~Hales, \emph{A proof of the Kepler conjecture},
Ann. of Math (2), \textbf{162} (2005), 1065--1185.

\bibitem[Hal06]{hales-2006}
T.C.~Hales, \emph{A proof of the Kepler conjecture (unabridged version)},
  Discrete Comput. Geom., \textbf{36} (2006), 1--265.

\bibitem[Hil01]{hilbert-1901}
D.~Hilbert, \emph{Mathematische {P}robleme}, Archiv Math. Physik \textbf{1}
  (1901), 44--63.

\bibitem[Hla49]{hlawka-1949}
E.~Hlawka, \emph{Ausf\"ullung und \"{U}berdeckung konvexer {K}\"orper durch
  konvexe {K}\"orper.}, Monatsh. Math. \textbf{53} (1949), 81--131.

\bibitem[Hsi01]{hsiang-2001}
W.Y.~Hsiang, \emph{{Least action principle of crystal formation of dense
  packing type and Kepler's conjecture.}}, World Scientific, Singapore, 2001.

\bibitem[Kep11]{kepler-1611}
J.~Kepler, \emph{Strena seu de nive sexangula}, Godefridum Tampach, Francofvrti
  Ad Moenum (Frankfurt), 1611, English translation: Hardy, {\em The
  six-cornered snowflake}, Clarendon Pr., Oxford, 1966.

\bibitem[Lag02]{lagarias-2002}
J.C.~Lagarias, \emph{Bounds for local density of sphere packings and the Kepler
  conjecture.}, Discrete Comput. Geom. \textbf{27} (2002), 165--193.

\bibitem[LG97]{lg-1997}
B.D.~Lubachevsky and R.L. Graham, \emph{Curved hexagonal packings of equal
  disks in a circle}, Discrete Comput. Geom. \textbf{18} (1997), 179--194.

\bibitem[Mel97]{melissen-1997}
J.B.M.~Melissen, \emph{Packing and covering with circles}, PhD thesis,
  University of Utrecht, 1997.

\bibitem[N\"O97]{no-1999}
K.J.~Nurmela and P.R.J.~\"Ostergard, \emph{More Optimal Packings of Equal Circles in a Square},
  Discrete Comput. Geom. \textbf{22} (1999), 439--457.

\bibitem[Sch93]{schneider-1993}
R.~Schneider, \emph{Convex bodies: the {B}runn-{M}inkowski theory.}, Cambridge
  University Press, Cambridge, 1993.

\bibitem[Sch02]{schuermann-2002}
A.~Sch\"urmann, \emph{On extremal finite packings}, Discrete Comput. Geom.
  \textbf{28} (2002), 389--403.

\bibitem[Spe04]{packomania}
E.~Specht, \emph{Packomania}, \url{http://www.packomania.com/}.

\bibitem[SMC+06]{sccscg-2006} P.G.~Szab\'{o}, M.Cs.~Mark\'{o}t, T.~Csendes, E.~Specht, L.G.~Casado, I.~Garc\'ia,
  \emph{New Approaches to Circle Packing in a Square},
  Springer-Verlag, Berlin, to appear.

\bibitem[Szp03]{szpiro-2003}
G.~Szpiro, \emph{Does the proof stack up?}, Nature \textbf{424} (2003), 12--13.

\bibitem[Zie97]{ziegler-1998}
G.M.~Ziegler, \emph{Lectures on polytopes}, Springer-Verlag, New York, 1997.

\bibitem[Zon99]{zong-1999}
C.~Zong, \emph{Sphere packings}, Springer-Verlag, New York, 1999.

\end{thebibliography}
\end{document}